\def\bbuildrel#1_#2^#3{\mathrel{\mathop{\kern 0pt#1}\limits_{#2}^{#3}}}

\def\NN{\mathbb N}

\def\QQ{\mathbb Q}

\def\Pot{\hbox{$\mathcal P$}}

\def\ms{\medskip}
\def\ss{\smallskip}
\def\noi{\noindent}

\def\w{\thinspace\hbox{\hsize 14pt \rightarrowfill}\thinspace}

\def\0{\hbox{$\emptyset$}}

\def\M{\hbox{\rm MGR}}

\def\bor{\mathbf B}
\def\s{\hbox{$\sigma$}}

\def\sub{\subseteq}

\def\M{\hbox{$\mathcal M$}}
\def\N{\hbox{$\mathcal N$}}
\def\F{\hbox{$\mathcal F$}}

\def\B{\mathscr{B}}

\def\E{\mathscr{E}}

\def\M{\mathscr{M}}

\def\U{\mathscr{U}}

\documentclass[a4paper,12pt]{amsart}

\usepackage{amssymb,enumerate,mathrsfs}
\usepackage[utf8]{inputenc} 

\theoremstyle{plain}

\newcommand{\cantor}{2^{\NN}}
\newcommand{\baire}{\NN^{\NN}}
\newtheorem{theorem}{Theorem}[section]

\newtheorem{lemma}[theorem]{Lemma}

\newtheorem{proposition}[theorem]{Proposition}

\numberwithin{equation}{section}

\begin{document}

\title{On countably perfectly meager
and countably perfectly null sets
}

\author{Tomasz Weiss and Piotr Zakrzewski}

\address{Institute of Mathematics, College of Science, Cardinal Stefan Wyszyński
University,  Dewajtis 5, 01-815 Warsaw, Poland}
\email{tomaszweiss@o2.pl}

\address{Institute of Mathematics, University of Warsaw,  Banacha 2, 02-097 Warsaw, Poland}
\email{piotrzak@mimuw.edu.pl}


\date{August 24, 2022}


\begin{abstract}

We study a strengthening of the notion of a universally meager set and its dual counterpart that strengthens  the notion of a universally null set.

We say that  a subset $A$ of   a perfect Polish space $X$ is countably perfectly meager (respectively, countably perfectly null) in $X$, if for every perfect Polish topology $\tau$ on $X$,  giving the original Borel structure of $X$, $A$ is covered by an  $F_\sigma$-set $F$ in $X$ with the original Polish topology such that $F$ is meager with respect to $\tau$ (respectively, for every finite, non-atomic, Borel measure $\mu$ on $X$, $A$ is covered by an $F_\sigma$-set $F$ in $X$  with $\mu(F)=0$).  

We prove that if $2^{\aleph_0}\leq \aleph_2$, then there exists a  universally meager set in $\cantor$ 
which is not countably perfectly meager in $\cantor$ (respectively, a universally null set in $\cantor$ 
which is not countably perfectly null in $\cantor$).

\end{abstract}

\maketitle

\section{Introduction}\label{sec:1} 

We continue the study of countably perfectly meager sets undertaken by Pol and Zakrzewski \cite{p-z}. We say (cf. \cite{p-z}) that a subset $A$ of  a perfect Polish
space $X$ is
{\it  countably perfectly meager in $X$} ($A\in \mathbf{PM}_\sigma$), if for  every sequence of perfect subsets $\{P_n: n \in \NN\}$ of $X$, there exists an $F_\sigma$-set $F$ in $X$ such
that $A \sub F$ and $F\cap P_n$ is meager in $P_n$ for each $n$.
Let us also recall that $A$    is {\it universally meager} ($A\in
\mathbf{UM}$), if  for every Borel
isomorphism $f$ between $X$ and 
any perfect Polish
space $Y$ the image of $A$ under $f$ is meager 
in $Y$  (see \cite{z-1}, \cite{z-2}, \cite{b-1}, \cite{b-2} and also \cite{g-1}, \cite{g-2}, \cite{g-3}, where this class was earlier
studied by Grzegorek and denoted by $\overline{\mathbf{AFC}}$).   By \cite[Theorem 7]{b-2} we have 
$\mathbf{PM}_\sigma\sub\mathbf{UM}$ and by \cite[Theorem 1.1]{p-z}, this inclusion is consistently proper, namely it holds if there exists a universally meager set of cardinality $2^{\aleph_0}$, in particular, if CH is true.

In this note we prove  (see Theorem \ref{PM_sigma neq UM}) that $\mathbf{PM}_\sigma\neq\mathbf{UM}$ follows also from the assumption that $2^{\aleph_0} = \aleph_2$. Whether it is consistent that
$\mathbf{PM}_\sigma=\mathbf{UM}$ remains an open problem (it is consistent that $\mathbf{UM}\subsetneq \mathbf{PM}$ but also that $\mathbf{UM}= \mathbf{PM}$ (see \cite{b-1}), where $\mathbf{PM}$ denotes the family of all perfectly meager subsets of $X$).

\ss

If $I$ is a \s-ideal of subsets of $X$, i.e., it is hereditary, closed under taking countable unions and contains all singletons, then by $I^*$ we  denote the \s-ideal on $X$ generated by the
closed subsets of $X$ which belong to $I$ (cf. \cite{s-z}).

\ss

If $\tau$ is a perfect Polish topology  on $X$  giving the original Borel structure of $X$, then by $\M(X,\tau)$ we denote the \s-ideal  of meager sets 
with respect to  $\tau$. Let us note that $\M^*(X,\tau)$ consists of such $A\sub X$ that there exists an $F_\sigma$-set $F$ in $X$ (with the original Polish topology) with $A\sub F$ and  $F\in\M(X,\tau)$. 
 By \cite[Theorem 2.1]{z-1}, $A$ is universally meager in $X$ if and only if $A$ belongs to the intersection of all \s-ideals of the form $\M(X,\tau)$,
whereas by \cite[Proposition 4.6]{p-z}, $A$ is countably perfectly meager in $X$ if and only if $A$ belongs to the intersection of all \s-ideals of the form $\M^*(X,\tau)$.

\ms

Universally meager sets may be seen  as a category counterpart of {\it universally null} sets in $X$. Namely, if for a finite, non-atomic, Borel measure $\mu$ is on $X$ (i.e., a countably additive measure $\mu: \bor(X)\w [0,+\infty)$ defined on the \s-algebra $\bor(X)$ of Borel subsets of $X$ and vanishing on singletons of $X$), we denote by $\N(X,\mu)$  the \s-ideal of $\mu$-null sets (i.e., sets of outer $\mu$-measure zero), then the collection  $\mathbf{UN}$ of universally null subsets of $X$ is the intersection of all \s-ideals of the form $\N(X,\mu)$.

\ss

The following definition of a measure analogue of countably perfectly meager sets was suggested by Taras Banakh. We say that $A$ is
{\it countably perfectly null in $X$} ($A\in \mathbf{PN}_\sigma$), if $A$ belongs to the intersection of all \s-ideals of the form $\N^*(X,\mu)$. In other words, $A\in \mathbf{PN}_\sigma$ if for every finite, non-atomic, Borel measure $\mu$ on $X$, $A$ is covered by an $F_\sigma$-set $F$ in $X$  with $\mu(F)=0$. Let us note that if $\lambda$ is  the standard probability product measure  on  the Cantor space $\cantor$, then $\N^*(\cantor,\lambda)$ is a well-known \s-ideal which is usually denoted by ${\E}$ (cf. \cite{b-j}).

 The name of the class $\mathbf{PN}_\sigma$ is further justified by the following observation.
 
 

 \begin{proposition}\label{characterisation}
 A set $A\sub X$ is countably perfectly null in $X$ if and only if 	for  every sequence of perfect subsets $\{P_n: n \in \NN\}$ of $X$ with associated probability non-atomic Borel measures $\mu_n$ on $P_n$, there exists an $F_\sigma$-set $F$ in $X$ such
 that $A \sub F$ and $\mu_n(F\cap P_n)=0$  for each $n$. 
 	\end{proposition}

 \begin{proof}
 If $A\in \mathbf{PN}_\sigma$ and for each $n$ we have	a perfect set $P_n$ together with the respective measure $\mu_n$ on $P_n$, then it is enough to cover $A$ by an $F_\sigma$-set $F$ with $\mu(F)=0$ for $\mu$ defined by
 $$
 \mu(B)=\sum_n \frac{1}{2^n}\mu_n(B\cap P_n)\quad\hbox{for}\quad B\in \bor(X).
 $$

 For the other direction, given a finite, non-atomic, Borel measure $\mu$ on $X$ let us note that the regularity of $\mu$ (cf. \cite[17.C]{ke}) implies the existence of (pairwise disjoint) perfect sets $\{P_n: n \in \NN\}$ of positive $\mu$-measure such that $\mu(X\setminus \bigcup_n P_n)=0$. Then it suffices to cover $A$ by an $F_\sigma$-set $F$ with $\mu(F\cap P_n)=0$ for each $n$.
 	
  \end{proof}

 Clearly, we have $\mathbf{PN}_\sigma \sub \mathbf{UN}$.  One easily observes that we also have 
 $\mathbf{PN}_\sigma \sub \mathbf{PM}_\sigma$.

\begin{proposition}\label{PNsubPM}
Every countably perfectly null subset of $X$ is countably perfectly meager.
\end{proposition}

\begin{proof}
Let us assume that $A\in \mathbf{PN}_\sigma$ and let $\{P_n: n \in \NN\}$ be a sequence of perfect subsets  of $X$. For each $n$ let $\mu_n$ be a Borel probability, non-atomic measure on $P_n$ which assigns positive values to all non-empty, relatively open subsets of $P_n$ (e.g., one may concentrate $\mu_n$ on a dense  in $P_n$ homeomorphic copy of  the irrationals). Let $F$ be 
an $F_\sigma$-set in $X$ such
that $A \sub F$ and $\mu_n(F\cap P_n)=0$  for each $n$ (cf. Proposition \ref{characterisation}). Clearly, $F\cap P_n$ is meager in $P_n$ for each $n$, so  $A\in \mathbf{PM}_\sigma$.
\end{proof}	
The inclusion $\mathbf{PN}_\sigma \sub \mathbf{PM}_\sigma$  is, at least consistently, proper. Indeed, if $A\sub \cantor$ is a Sierpiński set with respect to the measure $\lambda$, then $A\in\mathbf{PM}_\sigma$ in $\cantor$ (cf. \cite [Corollary 2.9 and Remark 2.11]{p-z}) but $A$ has positive outer measure $\lambda$.

An analogous argument shows the consistency of  $\mathbf{PN}_\sigma \neq \mathbf{UN}$. Namely, if $A\sub \cantor$ is a Luzin set in $\cantor$ (which exists e.g. under CH), then $A\in\mathbf{UN}$ but $A$ is not even  meager in $\cantor$.

In this note  we prove (see Theorem \ref{PN_sigma neq UN}) that the inequality $\mathbf{PN}_\sigma \neq \mathbf{UN}$ follows also from the assumptions that either there exists a universally null set in $\cantor$ of cardinality $2^{\aleph_0}$ (then we actually have that even $ \mathbf{UN} \setminus\mathbf{PM}_\sigma\neq\emptyset$; cf. Proposition \ref{PNsubPM})
or $2^{\aleph_0}= \aleph_2$. Whether it is consistent that
$\mathbf{PN}_\sigma=\mathbf{UN}$, remains an open problem. 

\ms

Section \ref{sec:2} is devoted to the proof of Theorem \ref{PM_sigma neq UM} stating that if $2^{\aleph_0}= \aleph_2$, then there is a universally meager set in $\cantor$ 
which is not countably perfectly meager in $\cantor$. 

In Section \ref{sec:3} we give some examples of countably perfectly null sets and prove Theorem \ref{PN_sigma neq UN} which shows the inequality $\mathbf{PN}_\sigma\neq\mathbf{UN}$ under the assumption that 
either there exists a universally null set in $\cantor$ of cardinality $2^{\aleph_0}$ (then we actually have that even $ \mathbf{UN} \setminus\mathbf{PM}_\sigma\neq\emptyset$, cf. Proposition \ref{PNsubPM})
or $2^{\aleph_0}= \aleph_2$.

In Section \ref{sec:4} we collect some remarks and open problems.

\section{Universally meager not countably perfectly meager sets}\label{sec:2}

Let us recall that the cardinal number $\mathfrak{b}$ is the minimal cardinality of a  subset of $\baire$ which is unbounded in the ordering  $\leq^*$  of eventual domination. Following \cite[Definition 2.8]{tsab}, by a {\sl $\mathfrak{b}$-scale} (in $\baire$) we mean
a  subset $B=\{f_\alpha:\alpha< \mathfrak{b}\}$ of 
$\baire$ with the following properties:
\begin{itemize}
	\item $f_\alpha:\NN\w\NN$ is strictly increasing,
	\item $\alpha<\beta<\mathfrak{b}$ implies $f_\alpha<^*f_\beta$,
	\item for every $f\in\baire$ there is $\alpha<\mathfrak{b}$ with $f_\alpha\nleq^*f$.	
\end{itemize}	

By identifying each $f_\alpha$ with the characteristic function of its  range (or just its range, respectively), we obtain a homeomorphic copy $A$ of $B$ in $\cantor$ (respectively, in $\Pot(\NN)$ with the Cantor set topology) which we also call a {\sl $\mathfrak{b}$-scale} in $\cantor$ (respectively, in $\Pot(\NN)$) (cf. \cite{tsab}). It is well-known and easy to see that $\mathfrak{b}$-scales can be  constructed in ZFC. They are also classical examples of sets which are both universally meager and universally null (cf. \cite{pl}).

Let us recall that given a  subset $A$ of a perfect Polish space $X$, by a {\sl $\gamma$–cover} of $A$ we mean a countable relatively open cover $\U$ of $A$ which is infinite and such that   
for each $x \in A$ the set $\{U \in\U : x \notin U\}$
is finite. We say that $A$ satisfies {\sl property $S_1(\Gamma,\Gamma)$} if  for every sequence $(\U_n : n\in\NN)$ of $\gamma$-covers of $A$ we can select
for each $n$ a set $V_n\in \U_n$
such that $\{V_n : n\in \NN\}$ is a $\gamma$-cover of $A$ (cf. \cite{j-m-s-s}, \cite{tsab}). It is well-known (and due to Hurewicz \cite{hur}) that property $S_1(\Gamma,\Gamma)$ implies the Hurewicz property (for a definition of the Hurewicz property see Section \ref{sec:3}).

 If $\mathfrak{b}=\omega_1$, then there exists a $\mathfrak{b}$-scale $A=\{a_\alpha:\alpha<\mathfrak{b}\}$ in $\Pot(\NN)$ with the additional property that $\alpha<\beta< \mathfrak{b}$ implies that $a_\beta \setminus a_\alpha$ is finite (see \cite[page 8]{tsab})  and by a theorem of Scheepers \cite{sch} (see also \cite[Theorem 123]{buk}), if  $A$ is such a $\mathfrak{b}$-scale in $\Pot(\NN)$, then $A\cup [\NN]^{<\aleph_0}$ has property $S_1(\Gamma,\Gamma)$. The following observation is an easy corollary of this result.  
Let us recall that if $\kappa$ is an infinite cardinal, then a set $A\sub X$ is $\kappa$-concentrated on a set $Q\sub X$, if $|A \setminus U| < \kappa$
for each open set $U$ in $X$ containing $Q$.

\begin{lemma}\label{S_1} Assume that $\mathfrak{b}=\omega_1$. Let $A=\{a_\alpha:\alpha<\mathfrak{b}\}$ be a $\mathfrak{b}$-scale  in $\Pot(\NN)$ with the additional property that $\alpha<\beta <\mathfrak{b}$ implies that $a_\beta \setminus a_\alpha$ is finite. 
	
	For each $n$ let 
	$\U_n=\{U^n_k: k\in\NN\}$ be an ascending (i.e., $U^n_{k}\sub U^n_{k+1}$) sequence  of open sets in $\Pot(\NN)$ with 
	$[\NN]^{<\aleph_0}\sub \bigcup_k U^n_k$ but 
	$[\NN]^{<\aleph_0}\sub U^n_k$ for no $k$. 
	Then we can select for each $n$ a set $V_n=U^n_{k_n}$
	such that   $\{V_n : n\in \NN\}$ is a $\gamma$-cover of $(A\cup [\NN]^{<\aleph_0})\setminus Y$ for a certain countable set $Y\sub A$.
	
\end{lemma}

\begin{proof}
	The set $A$ being a  $\mathfrak{b}$-scale in $\Pot(\NN)$, is 
	$\mathfrak{b}$-concentrated on $[\NN]^{<\aleph_0}$ (see \cite[Lemma 2.10]{tsab}). Consequently, since $\mathfrak{b}=\omega_1$, there is $\xi<\omega_1$ such that if we let $A'=\{a_\alpha:\xi<\alpha<\mathfrak{b}\}$, then for each $n$ we have $A'\cup [\NN]^{<\aleph_0}\sub \bigcup_k U^n_k$ and by the properties of  the sequence $\{U^n_k: k\in\NN\}$, $\{(A'\cup [\NN]^{<\aleph_0})\cap U^n_k: k\in\NN\}$ is a $\gamma$-cover of $A'\cup [\NN]^{<\aleph_0}$. Since at the same time $A'$ is still a $\mathfrak{b}$-scale  in $\Pot(\NN)$ with the additional property above, Scheepers's theorem gives the desired conclusion. 
\end{proof}

Let us recall that non($\M$) is the smallest cardinality of a non-meager subset of $\cantor$. It is well-knows that if $\tau$ is a perfect Polish topology  on a Polish space $X$, then non($\M$) is the smallest cardinality of a subset of $X$ not in $\M(X,\tau)$. We denote by $\QQ$  the copy of the rationals in $\cantor$
consisting of all eventually zero binary sequences.

\ss

Now we are ready to prove the main result of this section (cf. the proof of \cite[Theorem 4]{w}).

\begin{theorem}\label{PM_sigma neq UM}
	If $2^{\aleph_0}\leq \aleph_2$, then there is a universally meager set in $\cantor$ 
	which is not countably perfectly meager in $\cantor$.
	
\end{theorem}

\begin{proof}
	
	If 	$2^{\aleph_0}= \aleph_1$, then the result follows from \cite[Theorem 1.1]{p-z}, so from now on let us assume that $2^{\aleph_0}= \aleph_2$. 
	
	\ms 
	
	We shall split the argument into three cases.

	\ms
	
	
	{\bf Case (A):}  $\hbox{non}(\M)=\aleph_2$.
	
	Then, by a result of Grzegorek (see \cite[Theorem 1]{g-2}), there exists a universally meager set in $\cantor$ of cardinality $\aleph_2=2^{\aleph_0}$ and the result follows from \cite[Theorem 1.1]{p-z}.
	
	\ms

	{\bf Case (B):} $\mathfrak{b}=\aleph_2$.
	
	This case is  already covered by the previous one, since it is well-known that $\mathfrak{b}\leq \hbox{non}(\M)$.
	
	
	\ms

	{\bf Case (C):}  $\hbox{non}(\M)=\mathfrak{b}=\aleph_1$.
	
	\ms
	
	Let $C$ and $D$ be disjoint copies of the Cantor  set in $\cantor$ such that 
	\begin{enumerate}
		\item[(1)] the operation $+$ of addition is a homeomorphism between $C\times D$  and $C+D$ (cf. \cite{re-2}).
	\end{enumerate}
	
	Let us fix  a homeomorphism $h: \cantor \w C$. 
	
	Let $A=\{a_\alpha:\alpha<\mathfrak{b}\}$ be a $\mathfrak{b}$-scale  in $\Pot(\NN)$ with the additional property that $\alpha<\beta <\mathfrak{b}$ implies that $a_\beta \setminus a_\alpha$ is finite (cf. the paragraph preceding Lemma \ref{S_1}) and let us identify it with its homeomorphic copy  in $\cantor$.
	
	
	Let $X=A\cup \QQ$ and $\tilde{X}=h(X)$. Since $X$ is universally meager, so is $\tilde{X}$. 
	
	Let us fix a set $M\sub D$ of cardinality $\hbox{non}(\M)=\aleph_1$ such that 
	\begin{enumerate}
		\item[(2)]  $M$ is relatively non-meager in $D$. 	
	\end{enumerate}
	
	Since $|\tilde{X}|=\aleph_1$, we can fix a surjection  $m:\tilde{X}\w M$ onto $\tilde{X}$
	 and let $H=\{(x,m(x)):\ x\in\tilde{X}\}\sub C\times D$ be the graph of $m$.  Let us note that since $\tilde{X}$ is the injective continuous image of $H$ under the projection onto the first axis and
	$\tilde{X}$ is universally meager, so is $H$.
	
	Finally, let $Z=\{x+m(x):\ x\in\tilde{X}\}$. Clearly, $Z$ is universally meager as the image of $H$ under the homeomorphism $+$  between $C\times D$  and $C+D$ (cf. (1)). 
	\smallskip
	
	We shall show that 
	\begin{enumerate}
		\item[(3)] $Z$ is not a $\mathbf{PM}_\sigma$-set in $\cantor$
	\end{enumerate}	
	and this will end the proof of the theorem.
	
	To that end, let 
	$\tilde{\QQ}=
	h(\QQ)=\{q_n: n\in \NN\}$ and let us suppose, towards a contradiction, that there are closed  sets $F_n$ in $\cantor$ such that $Z\sub \bigcup_n F_n$ and $F_n$ is relatively nowhere dense in $q_k+D$ or equivalently,   $(q_k+F_n)\cap D$ is relatively nowhere dense in $D$ for each $n$ and $k$.
	
	Let 	$\{I_n:  n\in \NN\}$ be an enumeration with infinitely many repetitions 
	of  the elements of a countable basis $\B$
	of  $D$.
	
	Let us fix an arbitrary $i$ and let $F=F_i$.

	As the set $F$ is compact, for each $n$ we can define by induction on $k$  an ascending sequence $\{U^n_k: k\in\NN\}$ of open sets in $C$ with $\{q_i: i<k\}\sub U^n_k\cap\tilde{\QQ}\neq \tilde{\QQ}$ for every $k$ together with   a 
	sequence $\{D^n_k: k\in\NN\}$ of non-empty, relatively clopen sets  in $D$ such that
	\begin{enumerate}
		\item[(4)] $D^n_{k+1}\sub D^n_{k} \sub I_n$ and $cl_D((U^n_k+F)\cap D)\cap D^n_k=\emptyset$ for every $k$. 
	\end{enumerate}

	Now, since $\tilde{X}$ and $\tilde{\QQ}$ are the respective images of $X$ and $\QQ$ under the homeomorphism $h$, and  $\U_n=\{U^n_k: k\in\NN\}$ is an ascending sequence of open sets in $C$ with 
	$\tilde{\QQ}\sub \bigcup_k U^n_k$ but $\tilde{\QQ}\sub U^n_k$ for no $k$, Lemma \ref{S_1} enables us to
	select for each $n$ a set $V_n=U^n_{k_n}$
	such that
	\begin{enumerate}  
		\item[(5)]	$\{V_n : n\in \NN\}$ is a $\gamma$-cover of $\tilde{X}\setminus Y$
		for a certain countable set $Y\sub \tilde{X}$.
	\end{enumerate}	
	
	We will show that 
	\begin{enumerate}    
		\item[(6)]	  $((\tilde{X}\setminus Y)+F)\cap D$ is meager in $D$. 
	\end{enumerate}	
	
	To see this, for each $m$
	let  $K_m = \bigcap_{n\geq m} cl_D((V_n+F)\cap D)$ and let us note that $K_m$ is  a closed relatively nowhere dense subset of $D$. Indeed, any open set from $\B$ is of the form $I_n$ for some $n\geq m$  and  $I_n\not\subseteq K_m$ by (4).
	
	Moreover, we have
	$
	((\tilde{X}\setminus Y) + F)\cap D \sub \bigcup_m K_m.
	$
	Indeed, if $c\in  \tilde{X}\setminus Y$, then there is $m$ such that $c\in V_n$ for every $n\geq m$ (cf. (5)). Consequently, 
	$(c+F)\cap D\sub \bigcap_{n\geq m} ((V_n+F)\cap D)\sub K_m$, completing the proof of (6).
	
	\ms 
	
	Let us summarize: for each $i$ we have found a countable set $Y_i\sub \tilde{X}$ such that 
	$((\tilde{X}\setminus Y_i)+F_i)\cap D$ is meager in $D$.

	Consequently, letting $\tilde{Y} = \bigcup_i Y_i$ we get a countable subset of $C$ such that
	$((\tilde{X}\setminus\tilde{Y} )+\bigcup_n F_n)\cap D$ is meager in $D$.

	But since $Z\sub \bigcup_n F_n$, we conclude that
	\begin{enumerate}
		\item[(7)] $((\tilde{X}\setminus\tilde{Y} )+Z)\cap D$ is meager in $D$. 
	\end{enumerate}	
	
	On the other hand,  $M\setminus m(\tilde{Y})\sub (\tilde{X}\setminus\tilde{Y})+Z$. Indeed, if $m\in M\setminus m(\tilde{Y})$, then $m=m(x)$ for some $x\in \tilde{X}\setminus \tilde{Y}$ and then $m=(x+(x+m(x)))\in x+Z$. This implies that $((\tilde{X}\setminus\tilde{Y})+Z)\cap D$ is not meager in $D$
	(cf. (2)) contradicting (7) and thus completing the proof of (3).
	
\end{proof}	

Let us note that under CH we have $\hbox{non}(\M)=\mathfrak{b}=\aleph_1$ and  Case (C) of the proof above establishes the consistency of $\mathbf{PM}_\sigma\neq\mathbf{UM}$ in the way which avoids  the use of 
\cite[Theorem 1.1]{p-z}.

\section{Countably perfectly null sets}\label{sec:3} 

Let us recall that given a perfect Polish space $X$ a set $A\sub X$ has {\it the  Hurewicz property}, if for each sequence $\U_1, \U_2,\ldots$ of open  covers of $A$,
there are finite subfamilies $\F_n\sub \U_n$ such that $A\sub \bigcup_n \bigcap_{m\geq n}(\bigcup \F_m)$. If $A$ is a zero-dimensional subspace of $X$, then by a result of Hurewicz (cf. \cite{hur} and \cite{re-1}) this is equivalent  
 to the statement that every continuous image of $A$ in $\baire$ is bounded in the ordering  $\leq^*$  of eventual domination.

The smallest cardinality of a subset of $\cantor$ which is  nonmeasurable with respect to the standard probability product measure $\lambda$ on $\cantor$ is denoted by non($\N$). It is well-knows that if $\mu$ is a non-zero, finite, non-atomic, Borel measure on $X$, then non($\N$) is the smallest cardinality of a subset of $X$ not in $\N(X,\mu)$.

Let us also recall that by $\QQ$ we denote the copy of the rationals in $\cantor$
consisting of all eventually zero binary sequences.

\ss 

The following result provides  examples of universally null countably perfectly meager sets which are countably perfectly null as well.

\begin{proposition}\label{examples} The following collections of sets   are countably perfectly null in  the respective perfect Polish spaces:
	\begin{enumerate}
		\item universally null sets with the Hurewicz property in any perfect Polish space $X$,
		\item any sets of cardinality less than 
		$\min(\hbox{non}(\N), \mathfrak{b})$ in any perfect Polish space $X$,
		\item  $\gamma$-sets in any perfect Polish space $X$,
		\item $\mathfrak{b}$-scales in $\cantor$,
		
		\item  Hausdorff $(\omega_1,\omega_1^{*})$-gaps in $\Pot(\NN)$.
		
	\end{enumerate}
\end{proposition}

\begin{proof}
(1)	Let $A\sub X$ be a universally null set with the Hurewicz property and let $\mu$ be a non-zero, finite, non-atomic Borel measure on $X$. Since $A\in  \mathbf{UN}$, there is a $G_\delta$-set $G$ in $X$ such that $A\sub G$ and $\mu(G)=0$. Now, since $A$ has the Hurewicz property, there is an $F_\sigma$ set $F$ in $X$ such that $A\sub F\sub G$ (cf. \cite[Theorem 5.7]{j-m-s-s}). Consequently, $\mu(F)=0$ which shows that $A\in  \mathbf{PN}_\sigma$.

\ss

Statements (2) -- (4) can be derived from  (1) as follows.

\ss 
	
(2)	Sets of cardinality less than $\hbox{non}(\N)$ are universally null and  sets of cardinality less than $\mathfrak{b}$ have the Hurewicz property.

(3) $\gamma$-sets are  universally null (as they actually have Rothberger's property $C''$, cf. \cite{g-n}) and they have the Hurewicz property, by \cite[Theorem 2]{g-m}.
	
(4) Let us assume that $A$ is a $\mathfrak{b}$-scale in $\cantor$. Let $B=A\cup \QQ$. Then $B$ is a universally null set with the Hurewicz property (see e.g., \cite[Example 4.1 and Remark 4.2]{p-z}), so $B\in  \mathbf{PN}_\sigma$ in $\cantor$. Consequently, $A\in  \mathbf{PN}_\sigma$ in $\cantor$.
	
	(5). This may actually be established by a classical argument showing that the Hausdorff gap is universally null, which we sketch here for the sake of completeness. Following the proof of \cite[Lemma 20.5]{j-w}, let $\langle\langle a_\alpha:\alpha<\omega_1\rangle, \langle b_\alpha:\alpha<\omega_1\rangle \rangle$ be a Hausdorff gap,  $F_\alpha=\{c\in \Pot(\NN): a_\alpha \subseteq^* c\subseteq^*  b_\alpha \}$ for $\alpha<\omega_1$ and let $\mu$ be a non-zero, finite, non-atomic Borel measure on $\Pot(\NN)$. Then $F_\alpha$'s are $F_\sigma$-sets in $\Pot(\NN)$ and for a sufficiently large $\xi$ we have $\mu(F_\xi)=0$ (see \cite[the proof of Lemma 20.5]{j-w}). Letting 
	$$F=F_\xi \cup \{a_\alpha:\alpha<\xi\}\cup\{b_\alpha:\alpha<\xi \},$$
	we get an $F_\sigma$-set with $\{ a_\alpha:\alpha<\omega_1\}\cup \{ b_\alpha:\alpha<\omega_1\}\sub F$ and $\mu(F)=0$ which shows that $\{ a_\alpha:\alpha<\omega_1\}\cup \{ b_\alpha:\alpha<\omega_1\}\in  \mathbf{PN}_\sigma$ in $\Pot(\NN)$.

\end{proof}

The main result of this section is a measure counterpart of 
\cite[Theorem 1.1]{p-z} and Theorem \ref{PM_sigma neq UM}.

\begin{theorem}\label{PN_sigma neq UN}
	If either 
	\begin{enumerate}
		\item[(a)] there exists a universally null set in $\cantor$ of cardinality $2^{\aleph_0}$
	\end{enumerate}
or
\begin{enumerate}
	\item[(b)] $2^{\aleph_0}\leq \aleph_2$,
\end{enumerate}
then there is a universally null set in $\cantor$ 
which is not countably perfectly null in $\cantor$.

\end{theorem}

\begin{proof}
	
$(a)$	Let $T$ be a universally null set in $\cantor$ of cardinality $2^{\aleph_0}$.

By Proposition \ref{PNsubPM}, it suffices to show that 
there is also one which is not countably perfectly meager.

Let us recall that by \cite[Theorem 1.1]{p-z}, there exist a set $H\sub T\times \cantor$ intersecting each vertical section $\{t\}\times \cantor$, $t\in T$, in a singleton and a homeomorphic copy $E$ of $H$ in $\cantor$  which is not a $\mathbf{PM}_\sigma$-set  in $\cantor$. Now, since $T$ is 
universally null, so is $E$ as a preimage of $T$ under  
a continuous injective function.

\ss

$(b)$ If 	$2^{\aleph_0}= \aleph_1$, then  any Luzin set in $\cantor$ provides an example of a non-meager, universally null set.

From now on let us assume that $2^{\aleph_0}= \aleph_2$. 

\ms 

Following closely the scheme of proof of the Theorem \ref{PM_sigma neq UM}, we split the argument into three cases.

\ms


{\bf Case (A):}  $\hbox{non}(\N)=\aleph_2$.

Then, by a theorem of Grzegorek (see \cite{g}), there exists a universally null set in $\cantor$ of cardinality $\aleph_2=2^{\aleph_0}$ and the result follows from part (a).

\ms

{\bf Case (B):} $\mathfrak{b}=\aleph_2$.

In this case 
any  $\mathfrak{b}$-scale in $\cantor$ 
is  a universally null set of cardinality $\mathfrak{b}=2^{\aleph_0}$ and the result again follows from part (a).

\ms

{\bf Case (C):}  $\hbox{non}(\N)=\mathfrak{b}=\aleph_1$.

\ms

As  in the proof of Theorem \ref{PM_sigma neq UM}, we fix copies $C,\ D$ of the Cantor  set in $\cantor$ such that 
\begin{enumerate}
	\item[(1)] the operation $+$ of addition is a homeomorphism between $C\times D$  and $C+D$ (cf. \cite{re-2}),
\end{enumerate}
a homeomorphism $h: \cantor \w C$,  a  $\mathfrak{b}$-scale $X$ in $\cantor$ and we let $\tilde{X}=h(X)$. Since $X$ is universally null, so is $\tilde{X}$. 

We also fix  a homeomorphism $g: \cantor \w D$ and we define a Borel measure $\mu$ on $\cantor$ by letting
$$
\mu(B)=\lambda(g^{-1}(B\cap D)),\quad\hbox{for}\quad B\in \bor(\cantor).
$$

Then we fix a set $M\sub D$ of cardinality $\hbox{non}(\N)=\aleph_1$ with
\begin{enumerate}
	\item[(2)]  $\mu^*(M)>0$,	
\end{enumerate}
 we let $m:\tilde{X}\w M$ be a surjection onto $M$ and we put $H=\{(x,m(x)):\ x\in\tilde{X}\}$.  Since $\tilde{X}$ is the injective continuous image of $H$ under the projection onto the first axis and
$\tilde{X}$ is universally null, so is $H$.

Finally, let $Z=\{x+m(x):\ x\in\tilde{X}\}$. Clearly, $Z$ is universally null as the image of $H$ under the homeomorphism $+$  between $C\times D$  and $C+D$ (cf. (1)). 
	\smallskip

We shall show that on the other hand
\begin{enumerate}
	\item[(3)] $Z$ is not a $\mathbf{PN}_\sigma$-set in $\cantor$,
\end{enumerate}	
thus completing  the proof of the theorem.

To that end, let $\tilde{\QQ}=h(\QQ)=\{q_n: n\in \NN\}$ and let us suppose, towards a contradiction, that there are closed $\mu$-null sets $F_n$ in $\cantor$ such that $Z\sub \bigcup_n F_n$ and $\mu(q_k+F_n)=0$  for each $n$ and $k$ (cf. Proposition \ref{characterisation}.)

  Let us fix an arbitrary $\varepsilon>0$. 
 
 For each $n$, $F_n$ being compact and $\mu$-null, there is an open set $U_n$ in $C$ such that $\tilde{\QQ}\sub U_n$ and  
\begin{enumerate}
	\item[(4)] $\mu(U_n+F_n)<\frac{\varepsilon}{2^{n+1}}$.
\end{enumerate}

Now, $X$ being a  $\mathfrak{b}$-scale in $\cantor$, is 
$\mathfrak{b}$-concentrated on $\QQ$ (see \cite[Lemma 2.10]{tsab}). Consequently,  $\tilde{X}$ is 
$\mathfrak{b}$-concentrated on $\tilde{\QQ}$ which, taking into account that $\mathfrak{b}=\aleph_1$,  implies that for each $n$ there is a countable set $Y_n\sub \tilde{X}$ such that 
$\tilde{X}\setminus Y_n\sub U_n$. It follows (cf. (4)) that 
$\mu^*((\tilde{X}\setminus Y_n)+F_n)<\frac{\varepsilon}{2^{n+1}}$
 which implies that, letting $F=\bigcup_n F_n$ and $\tilde{Y}=\bigcup_n Y_n$, we have $\mu^*((\tilde{X}\setminus \tilde{Y})+F)<\varepsilon$. But since $Z\sub F$ and the choice of $\varepsilon$ was arbitrary, we conclude that
\begin{enumerate}
	\item[(5)]  $\mu^*((\tilde{X}\setminus \tilde{Y})+Z)=0$.
\end{enumerate}	

	On the other hand,   exactly as in the proof of Theorem \ref{PM_sigma neq UM}, we have $M\setminus m(\tilde{Y})\sub (\tilde{X}\setminus\tilde{Y})+Z$ 
	which, $\tilde{Y}$ being countable, implies that $\mu^*((\tilde{X}\setminus\tilde{Y})+Z)>0$ 
(cf. (2)), contradicting (5) and thus completing the proof of (3).

\end{proof}

\section{Remarks and open problems}\label{sec:4} 

The results of Sections \ref{sec:2} and \ref{sec:3} motivate the following questions. The first two are directly related to Theorems \ref{PM_sigma neq UM} and \ref{PN_sigma neq UN}, respectively.

\ss

\noi {\bf Problem 1.} Is $\mathbf{PM}_\sigma=\mathbf{UM}$ consistent?

\ss 
\noi {\bf Problem 2.} Is $\mathbf{PN}_\sigma=\mathbf{UN}$ consistent?

\ss

Let us note that  we consistently have $\mathbf{PM}_\sigma \sub \mathbf{UN}$ since  in the model obtained by adding $\aleph_2$ Cohen reals to a model of GCH we have $\mathbf{UM} \sub \mathbf{UN}$ (see Corazza \cite[Theorem 0.6(b)]{cor} and Miller \cite{mi}; by a theorem of Bartoszyński and Shelah, cf \cite[Theorem 3]{b-sh}, it is consistently true that even all perfectly meager sets are universally null). By the fact that $\mathbf{PN}_\sigma \sub\mathbf{PM}_\sigma$ (see Proposition \ref{PNsubPM}), the dual statement that  $\mathbf{PN}_\sigma \sub \mathbf{UM}$ is just true but the following question remains open.

\ss

\noi {\bf Problem 3.} Is $\mathbf{PN}_\sigma = \mathbf{PM}_\sigma$ consistent?

\ss

Finally, in view of the inclusion $\mathbf{PN}_\sigma \sub \mathbf{PM}_\sigma\cap \mathbf{UN}$ one may ask

\ss

\noi {\bf Problem 4.} Is $\mathbf{PN}_\sigma = \mathbf{PM}_\sigma\cap \mathbf{UN}$ true/consistent?



	


\bibliographystyle{amsplain}

\end{document}